\documentclass{amsart}

\usepackage{amssymb}
\usepackage{amsmath}
\usepackage{enumerate}

\def\s{\mathbb{S}}
\def\h{\mathbb{H}}
\def\r{\mathbb{R}}
\def\z{\mathbb{Z}}

\def\c{\mathbb{C}}

\newtheorem{remark}{Remark}
\newtheorem{theorem}{Theorem}
\newtheorem{proposition}{Proposition}

\newtheorem{example}{Example}

\numberwithin{equation}{section}

\begin{document}

\title[Minimal Lagrangian surfaces in $\s^2\times\s^2$]{Minimal Lagrangian surfaces in $\s^2\times\s^2$}

\author{Ildefonso Castro}
\address{Departamento de Matem\'{a}ticas \\
Universidad de Ja\'{e}n \\
23071 Ja\'{e}n, SPAIN}
\email{icastro@ujaen.es}
\thanks{Research partially supported by a MEC-Feder grant MTM2004-00109}

\author{Francisco Urbano}
\address{Departamento de Geometr\'{\i}a  y Topolog\'{\i}a \\
Universidad de Granada \\
18071 Granada, SPAIN}
\email{furbano@ugr.es}



\date{}

\begin{abstract}
We deal with the minimal Lagrangian surfaces of the Einstein-K\"ahler surface $\s^2\times\s^2$, 
studying local geometric properties 
and showing that they can be locally described as Gauss maps of minimal surfaces in $\s^3\subset \r^4$. 
We also discuss the second variation of the area and characterize the most relevant examples 
by their stability behaviour.
\end{abstract}

\maketitle

\section{Introduction}
The theory of minimal surfaces in three-dimensional Riemannian manifolds of constant sectional curvature is a classical topic in differential geometry which has been developed in large. Recently, many geometers are interested in minimal surfaces of $3$-dimensional manifolds of type $\Sigma\times\r$, where $\Sigma$ is a Riemann surface like the $2$-sphere $\s^2$ or the hyperbolic space $\h^2$.
As an illustrative example we refer to [MR] and references there in.

When the ambient space $M$ is a four-dimensional Riemannian manifold of constant sectional curvature, there are 
different approaches to the theory of minimal surfaces (see among others [B]) that is specially rich if $M$ is an Einstein-K\"ahler surface. 
In this case one can study their minimal surfaces according with their behaviour with respect to the K\"ahler structure, considering in this way important families of surfaces: complex surfaces (which are always minimal), Lagrangian surfaces, totally real surfaces, etc.\ (see for instance [W]).
The most relevant ambient space in this setting is the complex projective plane, whose minimal surfaces have been extensively studied.

Besides the complex projective plane, there is another only Hermitian symmetric space of compact type and complex dimension $2$: the Einstein-K\"ahler surface $\s^2\times\s^2$.
In this article we start the study of their minimal Lagrangian surfaces.  
In section 2, we give a brief introduction to the geometry of $\s^2\times\s^2$ as well as their Lagrangian surfaces, appearing a function $C$ on the surface (called {\em the associated Jacobian}) which is studied in depth and will play an important role along the paper. In section 3 we show two recipes for constructing Lagrangian surfaces of $\s^2\times\s^2$: as graphs of area-preserving diffeomorphisms of the sphere $\s^2$ (see Example 1) or as Gauss maps of certain surfaces of the Euclidean space $\r^4$ (see Example 2).

In section 4 we first classify the compact Lagrangian surfaces of $\s^2\times\s^2$ with non null parallel mean curvature vector (Theorem 1), by defining two holomorphic differentials on the surface that allow to prove that the associated Jacobian $C$ is an isoparametric function. We also study general properties of the minimal Lagrangian surfaces of $\s^2\times\s^2$, showing that the zeroes of the non negative function $1-4C^2$ are the zeroes of a holomorphic $2$-differential defined on the surface (see Proposition 5 for some consequences of this fact). Secondly we classify the minimal Lagrangian surfaces of $\s^2\times\s^2$ with constant Gauss curvature and provide rigidity results about the Gauss curvature of these surfaces (Theorem 2). 
The Gauss maps of orientable minimal surfaces of $\s^3\subset\r^4$ are weighty examples of minimal Lagrangian surfaces of $\s^2\times\s^2$ (see section 4.3). In Theorem 3 we prove a local converse of this fact, by establishing that any simply-connected minimal Lagrangian surface of $\s^2\times \s^2$ such that $1-4C^2$ has not zeroes is the Gauss map of a minimal surface of $\s^3\subset\r^4$.

A significant fact that also motivates the present paper is that the only example of extremal metric on a Klein bottle, i.e.\ a critical metric for the functional first eigenvalue of the Laplacian, recently discovered by D.\ Jakobson, N.\ Nadirashvili and I.\ Polterovich [JNP], can be described as the induced metric on the minimal Lagrangian Klein bottle embedded in $\s^2\times\s^2$ defined by $${\bf B}=\{((x,z),(y,w))\in \s^2\times\s^2 \,/\,2x=y, \ \Re(\sqrt{z}w)=\Im(\sqrt{z}w)\},$$
where we consider $\s^2=\{(x,z)\in\r\times\c\,/\,x^2+|z|^2=1\}.$ As a consequence, {\bf B} is a Hamiltonian stable minimal Lagrangian surface, i.e.\ a minimal Lagrangian surface stable for the area under Hamiltonian deformations of $\s^2\times\s^2$. 

In section 5 we study the second variation of the area functional for compact minimal Lagrangian surfaces of $\s^2\times\s^2$ proving the following uniqueness results:
\begin{quote}
{\em The totally geodesic Lagrangian sphere 
\[ {\bf M}_0=\{(x,-x)\in\s^2\times\s^2\,/\,x\in\s^2\} \]
is the unique stable minimal Lagrangian compact surface of $\s^2\times\s^2$.}
\end{quote}
\begin{quote}
{\em The totally geodesic Lagrangian sphere ${\bf M}_0$, the totally geodesic Lagrangian torus 
\[ {\bf T}=\{(x,y)\in\s^2\times\s^2\,/\,x_1=y_1=0\} \]
and the minimal Lagrangian Klein bottle {\bf B} are the unique Hamiltonian stable minimal Lagrangian compact surfaces of $\s^2\times\s^2$ with genus $g\leq 2$, when the surface is orientable, and with Euler characteristic $\chi\geq 0$, when the surface is non orientable.}
\end{quote}
\begin{quote}
{\em The index of an unstable minimal Lagrangian compact surface of $\s^2\times\s^2$ is at least $2$ and it is $2$ only for the totally geodesic Lagrangian torus {\bf T}.}
\end{quote}

\section{Lagrangian surfaces in $\s^2\times\s^2$}
Let $\s^2$ be the unit sphere in the Euclidean space $\r^3$ endowed with its standard Euclidean metric $\langle,\rangle$ and its structure of Riemann surface given by $J_xv=x\times v$, for any $v\in T_x\s^2$, $x\in\s^2$, where $\times$ stands for the vectorial product in $\r^3$. Its K\"{a}hler $2$-form is the area $2$-form $\omega_0$ defined by $\omega_0(v,w)=\langle J_x v,w\rangle=\det\{x,v,w\}$ for any $v,w\in T_x\s^2$.

We endow $\s^2\times\s^2$ with the product metric (also denoted by $\langle,\rangle$) and the product complex structure $J$ given by
\[
J_{(x,y)}(v)=(J_x v_1,J_y v_2)=(x \times v_1, y\times v_2),
\]
for any $v=(v_1,v_2)\in T_{(x,y)}(\s^2\times\s^2),\,(x,y)\in\s^2\times\s^2$, which becomes $\s^2\times\s^2$ in a K\"{a}hler surface. Its K\"{a}hler $2$-form is $\omega=\pi_1^*\omega_0+\pi_2^*\omega_0$ where $\pi_i$, $i=1,2$, are the projections of $\s^2\times\s^2$ onto $\s^2$.

The group of isometries of $\s^2\times\s^2$ is the subgroup of the orthogonal group $\hbox{O}(6)$ given by
\begin{equation}\label{isom}
\left\{ \left( \begin{array}{cc} A & 0 \\ 0 & B \end{array} \right) \,,\,\left( \begin{array}{cc} 0 & A \\ B & 0 \end{array} \right)\,/\, A,B \in \hbox{O}(3) \right\}.
\end{equation}
The subgroup of holomorphic (resp. subset of antiholomorphic) isometries of $\s^2\times\s^2$ is defined by the additional conditions $A,B\in \hbox{SO}(3)$ (resp.\ $\det A =\det B = -1 $). 
We point out that there are isometries of $\s^2\times\s^2$ neither holomorphic nor antiholomorphic. 

Considering $\s^2\times\s^2\subset\s^5(\sqrt 2)\subset\r^6$ and denoting by $\tilde{\sigma}$ the second fundamental form of $\s^2\times\s^2$ into $\r^6$, we have that
\begin{equation}\label{2ndff}
\tilde{\sigma}(v,w)=(-\langle v_1,w_1\rangle x,-\langle v_2,w_2\rangle y),
\end{equation}
where $v=(v_1,v_2), w=(w_1,w_2)\in T_{(x,y)}(\s^2\times\s^2)$, $(x,y)\in\s^2\times\s^2\subset\r^3\times\r^3$. 
It is clear that $\tilde{\sigma}$ satisfies 
\begin{equation}\label{2ndffJ}
\tilde{\sigma}(Jv,Jw)=\tilde{\sigma}(v,w),
\end{equation}
what implies that the mean curvature vector $\tilde{H}$ of $\s^2\times\s^2$ in $\r^6$ is given by
\begin{equation}\label{Htilde}
 2\tilde{H}_{(x,y)}=-(x,y), \ \forall (x,y)\in\s^2\times\s^2 .
\end{equation}
 In particular, $\s^2\times\s^2$ is a minimal hypersurface of $\s^5(\sqrt 2)$.  Using now the Gauss equation of $\s^2\times\s^2$ into $\r^6$, we obtain that $\s^2\times\s^2$ is an Einstein-K\"{a}hler surface of constant scalar curvature $4$. 

As an Hermitian symmetric space, we can identify $\s^2\times\s^2$ with the Grassmann manifold $G^+(2,4)$ of oriented $2$-planes in the Euclidean space $\r^4$ in the following way. Let $\Lambda^2\r^4=\{v\wedge w/v,w\in\r^4\}\equiv\r^6$ be the space of $2$-vectors in $\r^4$ endowed with the Euclidean metric $\langle\langle,\rangle\rangle$ given by 
\[
\langle\langle v\wedge w,v'\wedge w'\rangle\rangle=\langle v,v'\rangle\langle w,w'\rangle
-\langle v,w'\rangle\langle w,v'\rangle,
\]
for any $v,w,v',w'\in \r^4$.
We define the subspaces $\Lambda^2_{\pm}\r^4$ of $\Lambda^2\r^4$ generated by the unit vectors
\begin{eqnarray*}
E^1_{\pm}=\frac{\textstyle 1}{\textstyle \sqrt 2}(e_1\wedge e_2\pm e_3\wedge e_4),\\
E^2_{\pm}=\frac{\textstyle 1}{\textstyle \sqrt 2}(e_1\wedge e_3\pm e_4\wedge e_2), \\
E^3_{\pm}=\frac{\textstyle 1}{\textstyle \sqrt 2}(e_1\wedge e_4\pm e_2\wedge e_3),
\end{eqnarray*}
where $\{e_1,e_2,e_3,e_4\}$ is an oriented orthonormal frame of $\r^4$, and denote by $\s^2_{\pm}$ the unit spheres in the $3$-spaces $\Lambda^2_{\pm}\r^4$.

If $\{v_1,v_2\}$ is an oriented orthonormal frame of a plane $P\in G^+(2,4)$ and take $v_3$, $v_4$ in such a way that $\{v_1,v_2,v_3,v_4\}$ is an oriented orthonormal frame in $\r^4$, then the map 
\[
\begin{array}{c}
G^+(2,4)\longrightarrow \s^2_+\times\s^2_- \\
P \mapsto \left( \frac{\textstyle 1}{\textstyle \sqrt {2}}(v_1\wedge v_2+v_3\wedge v_4),
\frac{\textstyle 1}{\textstyle \sqrt {2}}(v_1\wedge v_2-v_3\wedge v_4)\right),
\end{array}
\]
defines a diffeomorphism. 

\vspace{0.3cm}

Let $\Phi=(\phi,\psi):\Sigma\rightarrow\s^2\times\s^2$ be an immersion of a surface $\Sigma$ and denote by $g=\phi^*\langle,\rangle +\psi^*\langle,\rangle$ the induced metric. The immersion $\Phi$ is said to be {\it Lagrangian} if $\Phi^*\omega=0$, i.e. $\phi^*\omega_0+\psi^*\omega_0=0$. This means that  
\[
0=\langle Jd\Phi_p(v),d\Phi_p(w)\rangle
= \langle Jd\phi_p(v),d\phi_p(w)\rangle+ \langle Jd\psi_p(v),d\psi_p(w)\rangle,
\]
for any $p\in\Sigma $ and $v,w\in T_p\Sigma$. 

If $\Phi=(\phi,\psi):\Sigma\rightarrow\s^2\times\s^2$ is an immersion of an oriented surface with area $2$-form $\omega_{\Sigma}$, 
we can define the Jacobians of $\phi$ and $\psi$ by
\[
\phi^*\omega_0=\hbox{Jac}\,(\phi)\,\omega_{\Sigma},\quad\psi^*\omega_0=\hbox{Jac}\,(\psi)\,\omega_{\Sigma}.
\]
Hence, when $\Sigma$ is oriented, {\em $\Phi$ is Lagrangian if and only if $\hbox{Jac}\,(\phi)=-\hbox{Jac}\,(\psi)$}. 
We will call the function 
\begin{equation}\label{defC}
C:=\hbox{Jac}\,(\phi)=-\hbox{Jac}\,(\psi)
\end{equation}
the {\em associated Jacobian} of the oriented Lagrangian surface $\Sigma$.
Moreover,  if $\Sigma$ is compact, we have that
\begin{equation}\label{defd}
\hbox{deg}\,(\phi)=-\hbox{deg}\,(\psi)=\frac{1}{4\pi}\int_{\Sigma}C\,\omega_{\Sigma}:=d.
\end{equation}
We will call this number $d$ the {\it degree} of the oriented compact Lagrangian surface $\Sigma$. 

In general (when $\Sigma$ is not necessarily orientable), the square $(\phi^*\omega_0)^2$ and $(\psi^*\omega_0)^2$ are well-defined and so it is also well-defined the function 
\[
\begin{array}{c}
C^2=\langle Jd\phi_p(e_1),d\phi_p(e_2)\rangle^2=\langle Jd\psi_p(e_1),d\psi_p(e_2)\rangle^2= \\ 
=|d\phi_p(e_1)\times d\phi_p(e_2)|^2=|d\psi_p(e_1)\times d\psi_p(e_2)|^2,
\end{array}
\]
where $\{e_1,e_2\}$ is an orthonormal basis of $(T_p\Sigma,g)$. 
If $\bar{R}$ denote the curvature operator of $\s^2\times\s^2$, from (\ref{2ndff}) it is easy to prove that $\bar{R}(e_1,e_2,e_2,e_1)=2C^2$. 
So the Gauss equation of $\Phi$ can be written as
\begin{equation}\label{EqG}
K=2C^2+2|H|^2-\frac{|\sigma|^2}{2},
\end{equation}
where $K$ is the Gauss curvature of $\Sigma$, $H$ the mean curvature of $\Phi$ and $\sigma$ the second fundamental form of $\Phi$. 
Thanks to the Lagrangian character of $\Phi $, the complex structure $J$ defines an isomorphism between the tangent bundle to $\Sigma $ and the normal bundle of $\Phi $. 
Therefore, $\{d\Phi_p(e_1),d\Phi_p(e_2),Jd\Phi_p(e_1),Jd\Phi_p(e_2)\}$ is an orthonormal frame at $T_{\Phi(p)}(\s^2\times\s^2)$, and using (\ref{2ndffJ}) and (\ref{Htilde}), we obtain that
\begin{equation}\label{Htilde2}
-2\Phi(p)=4\tilde{H}_{\Phi(p)}=2\sum_{i=1}^2 \tilde{\sigma}(d\Phi_p(e_i),d\Phi_p(e_i)).
\end{equation}
Since $\hat{\Phi}=(\phi,-\psi)$ is normal to $\s^2\times\s^2$ into $\s^5(\sqrt 2)$, we have that $\langle\tilde{H},\hat{\Phi}\rangle=0$.
If we use (\ref{2ndff}), (\ref{Htilde2}) and the fact that $\{e_1,e_2\}$ is an orthonormal frame in the above equality,  we prove that
\begin{equation}\label{ranks}
|d\phi_p(e_1)|^2+|d\phi_p(e_2)|^2=|d\psi_p(e_1)|^2+|d\psi_p(e_2)|^2=1.
\end{equation}
In particular, the maps $\phi,\psi:\Sigma\rightarrow\s^2$ satisfy that the ranks of $d\phi$ and $d\psi$ at any point of $\Sigma $ are always positive. Hence, if $C\equiv 0$ then the ranks of $d\phi$ and  $d\psi$ at any point must be necessarily $1$ and so both functions $\phi$ and $\psi$ define curves in $\s^2$. It is clear that the product of two spherical curves is a Lagrangian surface in $\s^2\times\s^2$ with null associated Jacobian.  
In conclusion, we get the following result.
\begin{proposition}\label{Prop:C=0}
Let $\phi:\Sigma\rightarrow\s^2\times\s^2$ be a Lagrangian immersion. Then $C\equiv 0$ if and only if $\Phi$ is locally the product immersion of two spherical curves
\[
\begin{array}{c}
I\times I' \longrightarrow \s^2\times\s^2 \\
(t,s) \mapsto (\alpha(t),\beta(s)).
\end{array}
\]
\end{proposition}
If we consider the product of two great circles of $\s^2$, we obtain a Lagrangian surface congruent to
\begin{equation}\label{defT}
{\bf T}=\{(x,y)\in\s^2\times\s^2\,/\,x_1=y_1=0\}\equiv\s^1\times\s^1.
\end{equation}
Of course, the associated Jacobian of $\bf T $ is null and its degree is zero. 
In addition, $\bf T$ is totally geodesic and flat.
 
\vspace{0.3cm}

We also get from (\ref{defC}) and (\ref{ranks}) that
\begin{eqnarray}\label{1-4C2}
\left(|d\phi_p(e_1)|^2-|d\phi_p(e_2)|^2\right)^2 +4\langle d\phi_p(e_1),d\phi_p(e_2)\rangle^2 = \\
=\left(|d\psi_p(e_1)|^2-|d\psi_p(e_2)|^2\right)^2 +4\langle d\psi_p(e_1),d\psi_p(e_2)\rangle^2=1-4C^2(p), \nonumber
\end{eqnarray}
which implies that $C^2(p)\leq 1/4$ and the equality holds if and only if 
$\phi$ and $\psi$ are conformal maps at $p\in\Sigma $.
If $C^2\equiv 1/4$, we have that $\phi$ and $\psi$ are conformal maps and (\ref{ranks}) implies that 
\begin{equation}\label{metrics}
\phi^*\langle,\rangle=\psi^*\langle,\rangle=g/2.
\end{equation} 
In particular, $\phi $ and $\psi $ are local diffeomorphisms. Thus, the immersion $\Phi $ can be locally reparametrized by
\[
\Phi=(i,F):U\subset\s^2\longrightarrow\s^2\times\s^2,
\]
where $i$ is the inclusion of $U$ into $\s^2$ and $F$ is a diffeomorphism from $U$ onto $F(U)$. 
But (\ref{metrics}) means that $F$ is the restriction to $U$ of an isometry $A$ of $\s^2$. Moreover, as $A^*\omega_0+\omega_0=0$ we have that $\det A=-1$. 
So our immersion is locally holomorphically congruent to $(i,-i):U\rightarrow\s^2\times\s^2$ by the holomorphic 
(see (\ref{isom})) isometry $\left( \begin{array}{cc} I & 0 \\ 0& -A \end{array}\right)$. We summarize this in the following result.
\begin{proposition}\label{Prop:C2=1/4}
 If $\Phi:\Sigma\rightarrow\s^2\times\s^2$ is a Lagrangian immersion then $C^2\leq 1/4$ and $C^2\equiv 1/4$ if and only if $\Phi(\Sigma)$ is congruent to an open subset of the Lagrangian surface ${\bf M}_0$ of $\s^2\times\s^2$ defined by
\begin{equation}\label{defM0}
{\bf M}_0=\{(x,-x)\in\s^2\times\s^2\,/\,x\in\s^2\}.
\end{equation}
\end{proposition}
It is easy to check that ${\bf M}_0$ is totally geodesic with $C\equiv 1/2$, $d=1$ and its Gauss curvature is constant $K\equiv 1/2$. 

In [CN] it was proved that {\em the totally geodesic Lagrangian surfaces of $\s^2\times\s^2$ are congruent to open subsets of $\,{\bf M}_0$ or $\bf T$.}

We study now the orientable compact Lagrangian surfaces embedded in $\s^2\times\s^2$.
\begin{proposition}\label{Prop:or+cpt+emb}
Let $\Phi:\Sigma\rightarrow\s^2\times\s^2$ be a Lagrangian immersion of an orientable compact surface $\Sigma$. If $\Phi$ is an embedding, then either the genus of $\Sigma$ is zero and the degree of $\Phi$ is $\pm 1$ or the genus of $\Sigma$ is one and the degree of $\Phi$ is zero.
\end{proposition}
\begin{proof}
In that follows, we use the notation of [W, Section 3]. 
We denote by $[\Sigma]\in H_2(\Sigma,\z)$ the fundamental homology class of $\Sigma$ and by $\Phi_*[\Sigma]^{\#}\in H^2(\s^2\times\s^2,\z)$ the Poincar\'{e} dual of $\Phi_*[\Sigma]\in H_2(\s^2\times\s^2,\z)$.  The self-intersection number of $\Phi$ is defined by
\[
\hbox{I}_{\Phi}=(\Phi_*[\Sigma]^{\#}\cup\Phi_*[\Sigma]^{\#})([\s^2\times\s^2]).
\]
Since $\Phi$ is an embedding, it is given by $\hbox{I}_{\Phi}=\chi(T^{\perp}\Sigma)$, where $\chi(T^{\perp}\Sigma)$ is the Euler number of the normal bundle. But 
\[
\begin{array}{c}
(\Phi_*[\Sigma]^{\#}\cup\Phi_*[\Sigma]^{\#})([\s^2\times\s^2])=(\Phi_*[\Sigma]^{\#})(\Phi_*[\Sigma]^{\#}\cap[\s^2\times\s^2])\\=(\Phi_*[\Sigma]^{\#})(\Phi_*[\Sigma]),
\end{array}
\]
and, using that $\Phi$ is a Lagrangian immersion, $\chi(T^{\perp}\Sigma)=-\chi(\Sigma)=2(g-1)$, where $g$ is the genus of $\Sigma$.
Hence, we obtain that
\begin{equation}\label{genus}
(\Phi_*[\Sigma]^{\#})(\Phi_*[\Sigma])=2(g-1).
\end{equation}

On the other hand, if $x$ is a point of $\s^2$, we know that $[\s^2\times\{x\}],[\{x\}\times\s^2]\in H_2(\s^2\times\s^2,\z)$ are generators of $H_2(\s^2\times\s^2,\z)\equiv\z^2$. Then it is well-known that their Poincar\'{e} duals satisfy:
\[
\begin{array}{c}
([\s^2\times\{x\}]^{\#})([\s^2\times\{x\}])=0,\quad ([\s^2\times\{x\}]^{\#})([\{x\}\times\s^2])=1, \\ 
([\{x\}\times\s^2]^{\#})([\s^2\times\{x\}])=1,\quad ([\{x\}\times\s^2]^{\#})([\{x\}\times\s^2])=0.
\end{array}
\]
As $\Phi_*[\Sigma]=\hbox{deg}\,(\phi)[\s^2\times\{x\}]+\hbox{deg}\,(\psi)[\{x\}\times\s^2]$, we obtain that
\[
([\Phi_*[\Sigma]^{\#})([\Phi_*[\Sigma])=2\,\hbox{deg}\,(\phi)\hbox{deg}\,(\psi).
\]
Using that $\hbox{deg}\,(\psi)=-\hbox{deg}\,(\phi)$, from (\ref{genus}) we finally get that $g-1=-\hbox{deg}\,(\phi)^2$, what proves the Proposition.
\end{proof}
\begin{remark} 
{\rm The totally geodesic Lagrangian surfaces ${\bf M}_0$ and $\bf T$ defined 
in (\ref{defM0}) and (\ref{defT}) show that the result proved in Proposition \ref{Prop:or+cpt+emb} is the best possible one.}
\end{remark}
In the following result, we show that ${\bf M}_0$ and $\bf T$ are the only examples of compact Lagrangian surfaces in $\s^2\times\s^2$ with constant associated Jacobian.
\begin{proposition}\label{Prop:C2cte}
There do not exist Lagrangian immersions of compact surfaces in $\s^2\times\s^2$ with constant $C^2\in(0,1/4)$.
\end{proposition}
\begin{proof}
Let $\Phi=(\phi,\psi):\Sigma\rightarrow\s^2\times\s^2$ be a Lagrangian immersion of a compact surface such that $C^2$ is constant, $0<C^2<1/4$. 
First, since $C^2$ does not have zeroes, the rank of $d\phi$ and $d\psi$ is $2$ at any point of $\Sigma$. 
So $\phi$ and $\psi$ are local diffeomorphisms and the compacity of $\Sigma$ implies that they are in fact diffeomorphisms. Therefore $\Sigma$ is a sphere, $\Phi$ must be an embedding and we can consider that $C$ (which is well-defined) is a constant $C\in (0,1/2)$. Then Proposition \ref{Prop:or+cpt+emb} says that the degree $d=1$ and so 
Area$(\Sigma)C=4\pi$ from (\ref{defd}). 

Next we prove that the Gauss curvature of $\Sigma $ is also constant, concretely $K=2C^2$:
Using the definition (\ref{defC}) of $C$, we have that
\[
2C(p)=\langle Jd\Phi_p(e_1),d\hat{\Phi}_p(e_2)\rangle,\,p\in\Sigma,
\]
where $\{e_1,e_2\}$ is an oriented orthonormal frame in $T_p\Sigma$ and $\hat{\Phi}=(\phi,-\psi)$. 
If $v\in T_p\Sigma$, then
\begin{equation}\label{vC}
2v(C)=\langle\sigma(v,e_2),Jd\hat{\Phi}_p(e_1)\rangle-\langle\sigma(v,e_1),Jd\hat{\Phi}_p(e_2)\rangle = 0.
\end{equation}
Using that the trilinear form $(u,v,w)\mapsto \langle \sigma (u,v), Jw \rangle$ is fully symmetric for a Lagrangian surface, we can write
\[
\sigma(e_1,e_1)=\lambda Je_1+ \eta Je_2,\,\,\sigma(e_1,e_2)=\eta Je_1+\mu Je_2,\,\,\sigma(e_2,e_2)=\mu Je_1+\nu Je_2,
\]
with $\lambda,\eta,\mu,\nu\in\r$.
Putting in (\ref{vC}) $v=e_1$ and $v=e_2$ respectively and using that $\{ e_1,e_2 \} $ is an orthonormal frame,
 we get that
\[
\begin{array}{c}
(\mu-\lambda)\,\langle d\phi_p(e_1),d\phi_p(e_2)\rangle + \eta (|d\phi_p(e_1)|^2-|d\phi_p(e_2)|^2)=0,\\ 
(\nu-\eta) \,\langle d\phi_p(e_1),d\phi_p(e_2)\rangle + \mu \,(|d\phi_p(e_1)|^2-|d\phi_p(e_2)|^2)=0.
\end{array}
\]
Taking into account that $1-4C^2(p)>0$,  from (\ref{1-4C2}) the last equations imply that $\mu(\mu-\lambda)=\eta(\nu-\eta)$. But this is equivalent to $|\sigma|^2=4|H|^2$ at any arbitrary point $p$. The Gauss equation (\ref{EqG}) now gives $K=2C^2$.
The Gauss-Bonnet theorem says that  $2\,$Area$(\Sigma)C^2=4\pi$, but we proved that Area$(\Sigma)C=4\pi$.
Thus we get $C=1/2$, which is a contradiction to the hypothesis.
\end{proof}

\section{Examples of Lagrangian surfaces.}
In this section, we emphasize two interesting ways of construction of Lagrangian surfaces in $\s^2 \times \s^2$.
\subsection{Graphs of area-preserving diffeomorphisms}
Let $F:U\subseteq\s^2\rightarrow\s^2$ be a smooth map defined on an open set $U$ of $\s^2$ and consider the graph of $F$,
\[
\begin{array}{c}
\Phi:U\longrightarrow\s^2\times\s^2\\
\Phi(x)=(x,F(x)).
\end{array}
\]
We have that $\Phi$ is Lagrangian if and only if $\omega_0 + F^* \omega_0 =0$.
In particular, $F$ must be a local diffeomorphism. 

When $U=\s^2$, $F$ must be a diffeomorphism and the last equation means that $-F$ preserves the area of $\s^2$. 
As a summary:

\begin{example}
If $F:\s^2\rightarrow\s^2$ is an area preserving diffeomorphism, the graph of $-F$:
\[
{\bf M}=\{(x,-F(x))\in \s^2\times\s^2\,/\,x\in\s^2\},
\]
is a Lagrangian surface of $\s^2\times\s^2$. 
In particular, the totally geodesic Lagrangian surface ${\bf M}_0$ defined in (\ref{defM0}) is the graph of the antipodal map $-I:\s^2\rightarrow\s^2$.
\end{example}

On the other hand, if $\Phi=(\phi,\psi):\Sigma\rightarrow\s^2\times\s^2$ is a Lagrangian immersion of a surface $\Sigma$ and $p\in\Sigma$ verifies $C^2(p)\not=0$, then there exists an open neighborhood $U$ around $p$ such that $\phi,\psi:U\rightarrow\s^2$ are diffeomorphisms onto their images, and so $\Phi(U)$ is the graph of $F=-\psi\circ\phi^{-1}$.

Globally, if $\Phi:\Sigma\rightarrow\s^2\times\s^2$ is a Lagrangian immersion of a connected compact surface $\Sigma$ such that $C^2$ has no zeroes, then $\Sigma$ is a sphere, $\Phi$ is an embedding and $\Phi(\Sigma)$ is the graph of the an area preserving diffeomorphism (with the opposite sign) of the sphere $\s^2$.
\begin{remark}
{\rm The result proved in Proposition \ref{Prop:C2cte} is not true if we do not assume the compacity of the surface, because we can construct examples of {\em non compact} Lagrangian surfaces in $\s^2\times\s^2$ whose associated Jacobian is constantly $\lambda$, for any $\lambda\in (0,1/2)$. In fact, we consider the map $F:\s^2-\{N,S\}\rightarrow\s^2$ defined by
\[
F(x,y,z)=(-e^{i\frac{\sqrt{1-4\lambda^2}}{\lambda}\tanh^{-1} z}(x+iy),z),
\]
where $N=(0,0,1)$ and $S=(0,0,-1)$. Then, it is straightforward to check that $F$ is an area preserving diffeomorphism from $\s^2-\{N,S\}$ onto itself, and so the graph of $-F$,
\[
\begin{array}{c}
\s^2-\{N,S\}\longrightarrow\s^2\times\s^2\\
(x,y,z)\mapsto \left((x,y,z),(e^{i\frac{\sqrt{1-4\lambda^2}}{\lambda}\tanh^{-1} z}(x+iy),-z)\right),
\end{array}
\]
is a Lagrangian embedding. It is an easy exercise to verify that the above Lagrangian graph has constant associated Jacobian $C=\lambda$.}
\end{remark}

\subsection{Gauss maps of certain surfaces of $\r^4$}
Let $\Psi:\Sigma\rightarrow\r^4$ be an immersion of an oriented surface $\Sigma$, and 
\[
\begin{array}{c}
\Phi:\Sigma\longrightarrow G^+(2,4)\\
\Phi (p)=d\Psi_p(T_p\Sigma),
\end{array}
\]
its Gauss map. If $\{e_1,e_2\}$ is an oriented orthonormal basis of $d\Psi_p(T_p\Sigma)$, taking into account the identification given in section 2, $\Phi$ can be written as
\[ \Phi=(\phi,\psi):\Sigma\rightarrow\s^2_+\times\s^2_-\subset\Lambda^2\r^4, \]
where
\[
\phi(p)=\frac{1}{\sqrt 2}(e_1\wedge e_2+e_3\wedge e_4)\equiv E^1_+(p),\, 
\psi(p)=\frac{1}{\sqrt 2}(e_1\wedge e_2-e_3\wedge e_4)\equiv E^1_-(p),
\]
being $\{e_1,e_2,e_3,e_4\}$ an oriented orthonormal frame of $\r^4$ at $\Psi(p)$.
Looking at $\Phi $ in $\Lambda^2\r^4\equiv \r^6$, we note that $\Phi(p)=\phi(p)+\psi(p)=\sqrt 2 (e_1\wedge e_2)$
while $\hat{\Phi}(p)=\phi(p)-\psi(p)=\sqrt 2 (e_3\wedge e_4)$.

For any vector $v\in T_p\Sigma$, it is easy to obtain that 
\begin{eqnarray}\label{difer}
\nonumber d\phi_p(v)=\left(\langle\hat{\sigma}(v,e_2),e_3\rangle+\langle\hat{\sigma}(v,e_1),e_4\rangle\right)E^2_+(p)\\
 +\left(\langle\hat{\sigma}(v,e_2),e_4\rangle-\langle\hat{\sigma}(v,e_1),e_3\rangle\right)E^3_+(p),\\
\nonumber d\psi_p(v)=\left(\langle\hat{\sigma}(v,e_2),e_3\rangle-\langle\hat{\sigma}(v,e_1),e_4\rangle\right)E^2_-(p)\\
\nonumber +\left(\langle\hat{\sigma}(v,e_2),e_4\rangle+\langle\hat{\sigma}(v,e_1),e_3\rangle\right)E^3_-(p),
\end{eqnarray}
where $\hat{\sigma}$ stands for the second fundamental form of the immersion $\Psi$. 
Hence, a point $p\in\Sigma$ verifies that $\hbox{dim}\, \hbox{Ker} \,d\Phi_p  > 0$ if and only if there exists a non null vector $v\in T_p\Sigma$ such that $\hat{\sigma}(v,w)=0, \forall w\in T_p\Sigma$, i.e. the index of relative nullity (see [CK]) of $\Psi $ at $p$ is positive.

Using that $JE^2_{\pm}(p)=E^3_{\pm}(p)$ and the Ricci equation of $\Psi$, from (\ref{difer}) we get that
\begin{eqnarray}\label{inner}
\langle d\Phi_p(v),d\Phi_p(w)\rangle = 2\left(2\langle\hat{\sigma}(v,w),\hat{H}\rangle-\hat{K}\langle v,w\rangle \right),\\
\langle Jd\Phi_p(v),d\Phi_p(w)\rangle = -2\hat{R}^{\perp}(v,w,e_3,e_4),\nonumber
\end{eqnarray}
$\forall\, v,w\in T_p\Sigma$, where $\hat{H}$, $\hat{K}$ and $\hat{R}^{\perp}$ are the mean curvature vector, the Gauss curvature and the normal curvature operator respectively associated to the immersion $\Psi$. 
Using (\ref{inner}) we arrive at the following conclusion:
\begin{example}
The Gauss map $\Phi:\Sigma\rightarrow\s^2_+\times\s^2_-$ of an immersion $\Psi:\Sigma\rightarrow\r^4$ of an oriented surface $\Sigma$ is a Lagrangian immersion if and only if $\Psi$ has flat normal connection and the index of relative nullity of $\Psi $ is zero. In particular, the Gauss map of any immersion $\Psi:\Sigma\rightarrow\s^3\subset\r^4$ of an oriented surface $\Sigma$ into the $3$-dimensional unit sphere is a Lagrangian immersion.
\end{example}
\begin{remark}\label{Re:Clifford}
{\rm If $\Psi $ is a minimal immersion of an oriented surface (with flat normal connection and negative Gauss curvature),
from (\ref{difer}) we get that both components $\phi $ and $\psi $ of its Gauss map $\Phi $ are conformal and (\ref{1-4C2}) says that $C^2\equiv 1/4$. In this case, $\Phi (\Sigma )$ must be congruent to an open set of ${\bf M}_0$ according to Proposition \ref{Prop:C2=1/4}.

If we consider the totally geodesic $\s^2\subset\s^3\subset\r^4$, then it is an exercise to check that its Gauss map is also the totally geodesic Lagrangian surface ${\bf M}_0$ defined in (\ref{defM0}). 

It is also easy to obtain that the Gauss map of the Clifford torus $\{ (z,w)\in \s^3\subset\c^2 \, / \, |z|=|w|=1/\sqrt 2 \}$ is given by $$ (z,w)\in \s^1(1/\sqrt 2) \times \s^1 (1/\sqrt2)\mapsto 2 \left( (0,-zw),(0,\bar{z}w) \right)\in \s^2 \times \s^2, $$ where $\s^2=\{(x,z)\in\r\times\c\,/\,x^2+|z|^2=1\}$.
We note it is a two fold covering of the totally geodesic Lagrangian torus $\bf T$ given in (\ref{defT}).}
\end{remark}

\section{Minimal Lagrangian surfaces}
Let $\Phi=(\phi,\psi):\Sigma\rightarrow\s^2\times\s^2$ be a Lagrangian immersion of a surface $\Sigma$. 
If $\{e_1,e_2\}$ is an orthonormal frame in $\Sigma$, then $\{e_1,e_2,Je_1,Je_2\}$ is an orthonormal frame in $\Phi^*T(\s^2\times\s^2)$ and, using (\ref{2ndffJ}) and (\ref{Htilde}), we deduce that the mean curvature 
vector $\bar{H}$ of $\Sigma$ into $\r^6$ is given by
\begin{equation}\label{mcR6}
\bar{H}=H+\tilde{H}=H-\frac{1}{2}\Phi. 
\end{equation}

We consider a local isothermal parameter $z=x+iy$ on $\Sigma$, in such a way that
\begin{eqnarray}\label{conformal}
\langle \Phi_z,\Phi_z\rangle=\langle \phi_z,\phi_z\rangle+\langle \psi_z,\psi_z\rangle=0,\\ 
|\Phi_z|^2=|\phi_z|^2+|\psi_z|^2=e^{2u}/2, \nonumber
\end{eqnarray}
where the derivatives respect to $z$ and $\bar{z}$ are given by
$\partial_z=\frac{1}{2}\left (\frac{\partial}{\partial x}-i\frac{\partial}{\partial y}\right )$,
$\partial_{\bar z}=\frac{1}{2}\left (\frac{\partial}{\partial x}+i\frac{\partial}{\partial y}\right )$.
Taking into  account that $\langle\Phi_{zz},J\Phi_{\bar z}\rangle=\langle\Phi_{z\bar z},J\Phi_{z}\rangle$, the Gauss equation of $\Phi$ and (\ref{mcR6}) imply that
\begin{equation}\label{Frenet}
\Phi_{z\bar z}=\frac{e^{2u}}{2}\left(H-\frac{1}{2}\Phi\right), 
\end{equation}
\[
\Phi_{zz}= 2u_z\Phi_z+\langle H,J\Phi_z\rangle J\Phi_{ z}+ 
2e^{-2u}\langle\Phi_{zz},J\Phi_z\rangle J\Phi_{\bar z}-\frac{1}{2}\langle\Phi_z, \hat{\Phi}_z\rangle\hat{\Phi},
\]
where $\hat{\Phi}=(\phi,-\psi)$. 

Also, from (\ref{ranks}) and (\ref{conformal}) we obtain that 
\begin{equation}\label{conformal2}
|\Phi_z|^2=2|\phi_z|^2=2|\psi_z|^2=e^{2u}/2.
\end{equation}
Finally, using (\ref{defC}) we can write 
\begin{equation}\label{flaC}
\langle J\Phi_{\bar z},\hat{\Phi}_{ z}\rangle=- ie^{2u}C
\end{equation} 
and deduce that
\begin{equation}\label{FrenetN}
\hat{\Phi}_z=2e^{-2u}\langle\Phi_z,\hat{\Phi}_z\rangle\Phi_{\bar{z}}-2iCJ\Phi_z.
\end{equation}
This yields
\begin{equation}\label{modulos}
e^{4u}(1-4C^2)=4|\langle\Phi_z,\hat{\Phi}_z\rangle|^2=16|\langle\phi_z,\phi_z\rangle|^2=16|\langle\psi_z,\psi_z\rangle|^2.
\end{equation}

\subsection{Lagrangian surfaces in $\s^2\times\s^2$ with parallel mean curvature vector}
We first classify the compact Lagrangian surfaces of $\s^2\times\s^2$ with non null parallel mean curvature vector.
\begin{theorem}\label{Th:cpt+Hparalelo}
Let $\Phi:\Sigma\rightarrow\s^2\times\s^2$ be a Lagrangian immersion with non null parallel mean curvature vector of a compact surface $\Sigma$. Then $\Phi$ is congruent to a finite covering of one of the embedded tori
\[
T_{a,b}=\{(x,y)\in\s^2\times\s^2\,/\,x_1=a,y_1=b\},
\]
with $a,b\in [0,1)$, $a^2+b^2>0$.
\end{theorem}
\begin{proof}
Since the mean curvature vector $H$ is a non null parallel vector field and the immersion is Lagrangian, $JH$ is also a non trivial parallel vector field on $\Sigma$. Hence the surface $\Sigma $ is flat. Using (\ref{EqG}), this  implies that $|\sigma|^2=4(|H|^2+C^2)$. Therefore, taking the two-fold oriented covering of $\Sigma$ if necessary, $\Sigma$ is a torus. 

Using (\ref{Frenet}), we can check that the $1$-differential $\Upsilon$ on $\Sigma$ defined by
\[
\Upsilon(z)=\langle H,J\Phi_z\rangle\otimes dz,
\]
is holomorphic. Since $\Upsilon$ does not vanish because $H$ is non null, we can normalize it by $\Upsilon=(1/2)(dz)$, i.e. $2\langle H,J\Phi_z\rangle=1$.
Derivating this equality and using again (\ref{Frenet}), we have that
\[
0=\langle H,J\Phi_z\rangle_z=\langle H,2u_zJ\Phi_z\rangle=u_z.
\]
This means that $u$ is constant, say $\mu\in \r$, 
and so $|H|^2=4 e^{-2\mu}\langle H,J\Phi_{\bar z} \rangle^2 = e^{-2\mu}$. 

We define now a $4$-differential on $\Sigma $ by
\[
\begin{array}{c}
\Xi(z)=\left (8\langle H,J\Phi_z\rangle\langle\Phi_{zz},J\Phi_z\rangle+\langle\Phi_z,\hat{\Phi}_z\rangle^2\right )\otimes(dz)^4\\=\left (4\langle\Phi_{zz},J\Phi_z\rangle+\langle\Phi_z,\hat{\Phi}_z\rangle^2\right )\otimes(dz)^4.
\end{array}
\]
Using (\ref{Frenet}), (\ref{FrenetN}) and the fact that $u$ is constant, it is easy to prove that $\Xi$ is also holomorphic. Since $\Sigma$ is a torus, $\Xi$ can be written as $\Xi=\lambda e^{i\theta}(dz)^4$, with $\lambda\geq 0$ and $\theta \in \r$.

Under these conditions, we are going to compute the gradient and the Laplacian of $C$.
From (\ref{flaC}), (\ref{Frenet}) and (\ref{FrenetN}), we obtain that
\[
iC_z=\frac{e^{-2\mu}}{2}\langle\Phi_z,\hat{\Phi}_z\rangle - 2e^{-4\mu}\langle\Phi_{zz},J\Phi_z\rangle\langle\Phi_{\bar z},\hat{\Phi}_{\bar z}\rangle .
\]
Since $|\nabla C|^2=4e^{-2\mu}|C_z|^2$,  using (\ref{conformal}) and (\ref{modulos}), we obtain that
\[
|\nabla C|^2= \frac{1-4C^2}{2}\left(e^{-2\mu}+2C^2 \right)
-8e^{-8\mu}\Re(\langle\Phi_{zz},J\Phi_z\rangle\langle\Phi_{\bar z},\hat{\Phi}_{\bar z}\rangle^2),
\]
where $\Re $ denotes real part.
But as $\Xi=\lambda e^{i\theta}(dz)^4$ we have that
\[
\lambda^2=e^{4\mu}+4e^{6\mu}C^2+\frac{e^{8\mu}}{16}\left( 1-4C^2 \right)^2+8\Re(\langle\Phi_{zz},J\Phi_z\rangle\langle\Phi_{\bar z},\hat{\Phi}_{\bar z}\rangle^2).
\]
The last two equations allow to get that
\begin{equation}\label{nablaC}
|\nabla \!C|^2\!=\!e^{-4\mu}(1-\lambda^2e^{-4\mu})\!+\!\frac{8e^{-2\mu}(1\!+\!4C^2)+(1\!-\!4C^2)(1\!+\!12C^2)}{16}.
\end{equation}
On the other hand, derivating $C_z$ with respect to $\bar z$ and using (\ref{Frenet}) and (\ref{FrenetN}), 
we obtain that
\[
C_{z\bar z}=-C\left(1+\frac{e^{2\mu}}{4}+e^{2\mu}C^2\right).
\]
Since $\Delta C =4e^{-2\mu} C_{z \bar z}$, we finally deduce that
\begin{equation}\label{DeltaC}
\Delta\,C=-C(1+4e^{-2\mu}+4C^2).
\end{equation}
Equations (\ref{nablaC}) and (\ref{DeltaC}) mean that the function $C$ is isoparametric.  

Then we follow an standard reasoning. 
We work on the open set $U$ where $\nabla C\not=0$. 
We are going to prove that $U=\emptyset$ and so $C$ must be constant. For this purpose, 
taking into account that $K=0$, the Bochner formula says that
\[
\frac{1}{2}\Delta|\nabla C|^2=\langle\nabla C,\nabla(\Delta C)\rangle+\sum_{i=1}^2|\nabla_{e_i}\nabla C|^2,
\]
where $\{e_1,e_2\}$ is any orthonormal frame on $U$. 
Using (\ref{nablaC}) and (\ref{DeltaC}), we can easily prove that
\[
\begin{array}{c}
\Delta|\nabla C|^2= C (1+4e^{-2\mu}-12C^2) \Delta C +(1+4e^{-2\mu}-36C^2)|\nabla C|^2,\\
\langle \nabla C,\nabla(\Delta C)\rangle=-(1+4e^{-2\mu}+12C^2)|\nabla C|^2.
\end{array}
\]
We now take on $U$ the orthonormal frame $\{e_1=\nabla C/|\nabla C|,e_2\}$. 
From (\ref{nablaC}) and (\ref{DeltaC}) it is not difficult to obtain that
\[
\sum_{i=1}^2|\nabla_{e_i}\nabla C|^2=C^2 \left(\frac{1}{2}+2e^{-2\mu}-6C^2\right)^2+C^2\left(\frac{3}{2}+6e^{-2\mu}-2C^2\right)^2.
\]
Using the last three expressions on the Bochner formula, after a long straightforward computation,
we finally arrive at
\[
0= 16 e^{8\mu} C^4 +(288 e^{6\mu}+72 e^{8\mu})C^2+48\lambda^2-48e^{4\mu}-24 e^{6\mu}-3e^{8\mu},
\]
which implies that $C$ is constant on each connected component of $U$. 
This is impossible and so $U=\emptyset$. Thus $\nabla\,C=0$ and $C$ must be constant. Looking at (\ref{DeltaC}) we deduce that $C\equiv 0$. 
Proposition \ref{Prop:C=0} says that $\Phi $ is locally the product of two spherical curves. In this family of surfaces, it is easy to check that $H$ is parallel and non null if and only if both curves have constant curvature non simultaneously zero, i.e. both curves are circles of $\s^2$. This finishes the proof.
\end{proof}

\subsection{Geometric properties}
Let $\Phi=(\phi,\psi):\Sigma\rightarrow\s^2\times\s^2$ be a minimal Lagrangian immersion of a surface $\Sigma$. 
The first equation of (\ref{Frenet}) implies that $\phi_{z\bar z}=-\frac{e^{2u}}{4}\phi$ and $\psi_{z\bar z}=-\frac{e^{2u}}{4}\psi$.  This means that $\phi,\psi:(\Sigma,g)\rightarrow (\s^2,\langle,\rangle)$ are harmonic maps. Thus the associated Hopf differential 
\[
\Theta(z)=\langle\phi_z,\phi_z\rangle\otimes(dz)^2=-\langle\psi_z,\psi_z\rangle \otimes (dz)^2=\frac{1}{2}\langle\Phi_z,\hat{\Phi}_z\rangle\otimes(dz)^2
\]
is holomorphic. Hence $\Theta$ vanishes identically or it has isolated zeroes. From (\ref{modulos}), we have that
\begin{equation}\label{modTheta}
16|\Theta|^2=e^{4u}(1-4C^2).
\end{equation}
 Therefore, either $C^2$ is constantly $1/4$ and, 
according to Proposition \ref{Prop:C2=1/4}, our surface is congruent to an open subset of ${\bf M}_0$, 
or the points where $C^2=1/4$ are isolated. 

On the other hand, using some properties of harmonic maps studied by Schoen and Yau in [SY] (see also Chapter 5 of [J]) we have that, if $\Theta$ is not identically zero, the functions $(1+2C)/4$ and $(1-2C)/4$, outside the isolated points where they vanish, satisfy (see Lemma 5.2.1 in [J] or formulae (16) and (18) in [SY]):
\begin{equation}\label{Deltalogs}
\Delta\log \frac{1+2C}{4}=2K-2C,\quad \Delta\log \frac{1-2C}{4}=2K+2C.
\end{equation}
Taking into account all this information, we obtain the following result.
\begin{proposition}\label{Prop:min+cpt+or}
Let $\Phi:\Sigma\rightarrow\s^2\times\s^2$ be a minimal Lagrangian immersion of a surface $\Sigma$. Then:
\begin{enumerate}
\item Either $C^2\equiv1/4$ and $\Phi(\Sigma)$ is an open set of the totally geodesic Lagrangian surface ${\bf M}_0$ 
 or the points of $\Sigma$ where $C^2=1/4$ are isolated.
\item If $\Sigma$ is a sphere, then $\Phi(\Sigma)$ is congruent to ${\bf M}_0$. 
In particular, a real projective plane cannot be immersed in $\s^2\times\s^2$ as a minimal Lagrangian surface.
\item If $\Sigma$ is a compact orientable surface of genus $g\geq 1$, then the degree $d$ of $\Phi$ 
 satisfies
\[
1-g-d=-N^-/4,\quad 1-g+d=-N^+/4,
\]
being $N^-$ (respectively $N^+$) the sum of all orders for all zeroes of $1+2C$ (respectively $1-2C$). 
In particular, if $\Sigma$ is a torus, then $C^2<1/4$ and the degree of $\Phi$ is zero, i.e., $\int_{\Sigma}C\,dv_g=0$. 
\end{enumerate}
\end{proposition}
\begin{proof}
Part (1) is clear. 
If $\Sigma$ is a sphere, the Riemann-Roch theorem says that $\Theta\equiv 0$, which means $C^2\equiv1/4$ and part (2) follows. Part (3) is proved integrating (\ref{Deltalogs}) and using Lemma 5.2.3 in [J] or formulae (6) and (7) in [SY].
\end{proof}

Next we prove some general properties of the minimal Lagrangian surfaces of $\s^2\times\s^2$ related with their Gauss curvatures.
\begin{theorem}\label{Th:min+K}
Let $\Phi:\Sigma\rightarrow\s^2\times\s^2$ be a minimal Lagrangian immersion of a surface $\Sigma$. 
\begin{enumerate}
\item If the Gauss curvature $K$ of $\Sigma$ is constant, then $\Phi(\Sigma)$ is congruent to some open subset of the totally geodesic Lagrangian surfaces ${\bf M}_0$ or $\bf T$.
\item If $\Sigma$ is complete and the Gauss curvature $K$ is non negative, then $\Phi(\Sigma)$ is congruent to the sphere ${\bf M}_0$ or to the torus $\bf T$. 
\item If $\Sigma$ is complete, the Gauss curvature $K$ is non positive and $1-4C^2\geq\epsilon>0$ for some constant $\epsilon$, then $\Phi(\Sigma)$ is congruent to $\bf T$.
\end{enumerate}
\end{theorem}
\begin{proof}
From (\ref{flaC}) and (\ref{modulos}) we get that 
\begin{equation}\label{nablaDeltaC}
|\nabla C|^2=\frac{(1-4C^2)(2C^2-K)}{2},\quad \Delta C=-C(1+4C^2-4K).
\end{equation}
Suppose now that $K\equiv a$, with $a\in\r$. Then (\ref{nablaDeltaC}) says that the function $C$ is isoparametric. 
We follow a similar reasoning to the used in the proof of Theorem \ref{Th:cpt+Hparalelo}. 
We work on the open set $U$ where $\nabla C\not=0$. 
We are going to prove that $U=\emptyset$ and so $C$ must be constant.
Using that $K\equiv a$, the Bochner formula gives
\[
\frac{1}{2}\Delta|\nabla C|^2=a |\nabla C|^2+\langle\nabla C,\nabla(\Delta C)\rangle+\sum_{i=1}^2|\nabla_{e_i}\nabla C|^2,
\]
where $\{e_1,e_2\}$ is any orthonormal frame on $U$. Using (\ref{nablaDeltaC}), it is easy to prove that
\[
\begin{array}{c}
\Delta|\nabla C|^2=2C(1+2a-8C^2)\Delta C+2(1+2a-24C^2)|\nabla C|^2,\\
\langle \nabla C,\nabla(\Delta C)\rangle=-(1-4a+12C^2)|\nabla C|^2.
\end{array}
\]
We take on $U$ the orthonormal frame $\{e_1=\nabla C/|\nabla C|,e_2\}$. From (\ref{nablaDeltaC}) we obtain that
\[
\sum_{i=1}^2|\nabla_{e_i}\nabla C|^2=C^2(1+2a-8C^2)^2+C^2(2-2a-4C^2)^2.
\]
Using the last three expressions on the Bochner formula, we finally arrive at
\[
0=(1-4C^2) \left( (3a^2/2-a)+(9a-4)C^2 \right),
\]
which implies that $C$ is constant on each connected component of $U$. 
This is impossible and so $U=\emptyset$. Hence $\nabla C=0$ and so $C$ is constant. 
But using once again (\ref{nablaDeltaC}) we obtain that $C=0$ or $C^2=1/4$.
In the first case, Proposition \ref{Prop:C=0} says that $\Phi $ is locally the product of two spherical curves. In this family of surfaces, it is easy to check that $H$ is null if and only if both curves have zero curvature, i.e. both curves are great circles of $\s^2$ and this leads to $\bf T$.
In the second case, Proposition \ref{Prop:C2=1/4} leads to ${\bf M}_0$. This finishes the proof of part (1).

In order to prove parts (2) and (3), we first compute the Laplacian of the non negative function $1-4C^2$. In fact, 
from (\ref{nablaDeltaC}), it is clear that
\begin{equation}\label{Delta1}
\Delta (1-4C^2)=4K(1-4C^2)+16C^2|\sigma|^2.
\end{equation} 
If $K\geq 0$, then (\ref{Delta1}) says that $1-4C^2$ is a subharmonic function that satisfies $1-4C^2\leq 1$.
Since $\Sigma $ is a complete non negative curved surface, the maximum principle implies that $1-4C^2$ is constant. 
From (\ref{Delta1}) we have that $K=C \equiv 0 $ or $\sigma \equiv 0$.
This finishes the proof of part (2).

To prove part (3), we make use of (\ref{nablaDeltaC}) in order to obtain
\begin{equation}\label{Delta2}
\Delta\log\,(1-4C^2)=4K.
\end{equation}
Then (\ref{Delta2}) implies that $g_0=(1-4C^2)^{1/2}g$ is a flat and complete metric and hence $\Sigma$ is parabolic. 
If $\Delta_0$ is the Laplacian of $g_0$, then we deduce from (\ref{Delta2}) that
\[
\Delta_0\log(1-4C^2)=\frac{4K}{(1-4C^2)^{1/2}}\leq 0.
\]
We have that $\log (1-4C^2)$ is a superharmonic function and is indeed bounded away from zero if $1-4C^2 \geq \epsilon > 0$. So it must be constant and so $K\equiv 0 $ what finishes the proof.
\end{proof}

\subsection{Minimal Lagrangian surfaces as Gauss maps}
We emphasize important examples of minimal Lagrangian surfaces of $\s^2\times\s^2$. 
Let $\Psi:\Sigma\rightarrow\r^4$ be an immersion of an oriented surface $\Sigma$ and $\Phi=(\phi,\psi):\Sigma\rightarrow\s^2_+\times\s^2_-$ its Gauss map (see section 3.2). 
It is well-known (see [RV]) that $\Phi$ is harmonic if and only if $\Psi$ has parallel mean curvature vector. 
From (\ref{inner}), $\Phi$ is a conformal map if and only if $A_{\hat{H}}=|\hat{H}|^2 I$, where $A$ is the Weingarten endomorphism of the immersion $\Psi$ and $I$ denotes the identity map. 
Thus, $\Phi$ is minimal if and only if $\Psi$ has parallel mean curvature vector and $A_{\hat{H}}=|\hat{H}|^2 I$. 
If $\Psi $ is minimal, then $\Phi(\Sigma )$ must be an open set in ${\bf M}_0$ (see Remark \ref{Re:Clifford}).
If $\hat H $ is non null,
we deduce that $\Psi $ is a minimal immersion in a 3-dimensional sphere of radius $1/|\hat{H}|$. Up to translations and scale, we can consider that $\Psi $ lies in the 3-dimensional sphere of radius $1$ centred at $0$. In conclusion:
\begin{quote}
{\em The Gauss map $\Phi=(\phi,\psi):\Sigma\rightarrow\s^2_+\times\s^2_-$ of a minimal immersion $\Psi:\Sigma\rightarrow\s^3$ of an oriented surface $\Sigma$ in the $3$-dimensional unit sphere  
is a minimal Lagrangian immersion in $\s^2\times\s^2$.}
\end{quote}

We focus our attention in this last case. 
Let $\Psi:\Sigma\rightarrow\s^3$ be a minimal immersion of an oriented surface in $\s^3$ and let $\hat{\sigma}$ be now
its second fundamental form. The Gauss equation of $\Psi$ is written as $\hat{K}=1-|\hat{\sigma}|^2/2$.
We obtain from (\ref{inner})  that the induced metrics $g$ and $\hat g$ on $\Sigma$ by the immersions 
$\Phi$ and $\Psi$ respectively are conformal.  Concretely:
\begin{equation}\label{gGauss}
g=(2+|\hat{\sigma}|^2)\hat{g}.
\end{equation}
Moreover, using (\ref{difer}), (\ref{gGauss}) and the Gauss equation of $\Psi $, it is not difficult to check that the associated Jacobian $C$ of the minimal Lagrangian immersion $\Phi$ coming from the Gauss map of $\Psi$ (which is defined by $C(p)=\langle Jd\phi_p(v),d\phi_p(w)\rangle$, where $\{v,w\}$ is an oriented orthonormal frame in $(\Sigma,g)$) is given by
\begin{equation}\label{CGauss}
C=\frac{2-|\hat{\sigma}|^2}{2(2+|\hat{\sigma}|^2)}=\frac{\hat{K}}{2+|\hat{\sigma}|^2}\,.
\end{equation}
Hence, from (\ref{CGauss}) we deduce that $-1/2<C\leq 1/2$ and the points where $C=1/2$ correspond to the isolated zeroes of $\hat{\sigma}$.
In addition, if $\Sigma$ is compact, using (\ref{defd}), (\ref{gGauss}), (\ref{CGauss}) and the Gauss-Bonnet theorem, 
its degree $d$ is given by
\[
d=\frac{1}{4\pi}\int_{\Sigma}C\,dv_{g}=\frac{1}{4\pi}\int_{\Sigma}\hat{K}\,dv_{\hat{g}}=1-g,
\]
where $g$ is the genus of $\Sigma$. In particular (see Proposition \ref{Prop:min+cpt+or}), $N^-=0$ and $N^+=8(g-1)$.
 
On the other hand, another ''Gauss map'' $N:\Sigma \rightarrow \s^3$ is defined pointwise for $\Psi:\Sigma\rightarrow\s^3\subset \r^4$ as the image of the unit normal in $\s^3$ translated to the origin in $\r^4$. The image $N(\Sigma )$ is called a polar variety in [L] and it is a minimal surface with singularities occurring at the points where $\hat{K}= 1$. If we consider the Gauss map of $N$, we find a minimal Lagrangian immersion in $\s^2\times \s^2$ which induces the same metric that $\Psi $ but changes the sign of $C$. One can check that it is exactly the Gauss map $\Phi $ of $\Psi $ with the opposite sign.
Moreover, we can choose $N$ in such a way that the pair $\{\Psi,N\}$ is oriented in a compatible way so that (see section 3.2)  $\hat{\Phi}=(\phi,-\psi)=\sqrt 2(\Psi\wedge N)$. The immersion 
\[
\hat{\Phi}:\Sigma\rightarrow\s^2_+\times\s^2_-\subset\s^5(\sqrt 2)
\]
defines (see [L]) a minimal immersion in $\s^5(\sqrt 2)$ which is known as the bipolar of $\Psi$. 
Thus we deduce that the Gauss map $\Phi=(\phi,\psi)$ and the bipolar $\hat{\Phi}=(\phi,-\psi)$ of the immersion $\Psi$ are congruent immersions
since $\hat{\Phi}= {\bf I} \Phi$ but we point out that the isometry
$\bf I=\left( \begin{array}{cc} I & 0 \\ 0& -I \end{array}\right)$ is neither holomorphic nor antiholomorphic (see section 1).

Thanks to the work of [L] and using this procedure we can assert:
\begin{quote}
{\em
Every compact Riemann surface of arbitrary genus can be immersed in $\s^2\times\s^2$ as a minimal Lagrangian surface.}
\end{quote}

The next result shows that, outside a set of isolated points, any minimal Lagrangian surface in $\s^2\times\s^2$ is locally the Gauss map of a minimal oriented surface in $\s^3$.

\begin{theorem}\label{Th:min}
Let $\Phi:\Sigma\rightarrow\s^2\times\s^2$ be a minimal Lagrangian immersion of a simply-connected surface $\Sigma$ with $C^2<1/4$. Then $\Phi$ is congruent to the Gauss map of a minimal immersion $\Psi:\Sigma\rightarrow\s^3\subset\r^4$ in the $3$-dimensional unit sphere of $\r^4$.
\end{theorem}
\begin{proof}
Since $\Phi $ is minimal and $C^2<1/4$, from (\ref{modTheta}) we know that $\Theta(z)=\frac{\langle\Phi_z,\hat{\Phi}_z\rangle}{2}\otimes(dz)^2$ is a holomorphic $2$-differential without zeroes. So, up to a change of complex coordinate on $\Sigma$ if necessary, we can normalize it on $\Sigma$ by 
$\Theta(z)=e^{i\theta}(dz)^2$, $\theta \in \r$. This implies that
\begin{equation}\label{F}
\langle\Phi_z,\hat{\Phi}_z\rangle=2e^{i\theta},
\end{equation}
and (\ref{modulos}) gives
\begin{equation}\label{matrica}
1-4C^2=16e^{-4u}.
\end{equation}
Then (\ref{Frenet}), (\ref{flaC}) and (\ref{F}) lead to
\[
ie^{2u}C_z=4e^{-i\theta}e^{-2u}\langle\Phi_{zz},J\Phi_z\rangle,
\]
so that we obtain that
\[
e^{8u}|C_z|^2=16|\langle\Phi_{zz},J\Phi_z\rangle|^2=e^{6u}|\sigma|^2.
\]
Using this in the Gauss equation (\ref{EqG}), we deduce that the Gauss curvature $K$ of the surface is given by $K=2C^2-e^{2u}|C_z|^2/2$. Since $4u_{z\bar z}=-e^{2u}K$, from (\ref{Frenet}) and (\ref{FrenetN}) we finally can reach that the Gauss and Codazzi equations of the Lagrangian immersion $\Phi$ are equivalent to
\begin{equation}\label{GCeqns}
2u_{z\bar z}+e^{2u}C^2-\frac{e^{4u}|C_z|^2}{4}=0,\quad\quad 1-4C^2=16e^{-4u}.
\end{equation}
Taking into account that $-1<2C<1$, we can now define the function $v$ on $\Sigma$ by $\tanh 2v=2C$. 
Using (\ref{matrica}) we have that $e^{2u}=4\cosh 2v$ and
it is easy to check that equations (\ref{GCeqns}) are equivalent to the sinh-Gordon equation
\begin{equation}\label{sinhGordon}
v_{z\bar z}+\frac{\sinh 2v}{2}=0.
\end{equation}
It is well known that for any solution $v$ of (\ref{sinhGordon}) 
there exists a one-parameter family $\Psi^\vartheta:(\Sigma,e^{2v}|dz|^2)\rightarrow\s^3$, $\vartheta \in \r$,
of minimal isometric immersions of our simply connected surface, where its associated holomorphic $2$-differential $\Xi^\vartheta(z)=\langle\Psi^\vartheta_z,N^\vartheta_z\rangle\otimes(dz)^2$, with $N^\vartheta$ the unit normal to $\Psi^\vartheta$ such that $\{ \Psi^\vartheta_x, \Psi^\vartheta_y, \Psi^\vartheta , N^\vartheta \}$ is an oriented
frame, has been normalized by $\Xi=e^{i\vartheta}(dz)^2/2$.

We study now the Gauss map $\Phi^\vartheta $ of $\Psi^\vartheta$.  
The Gauss equation of $\Psi^\vartheta$ gives that the norm of its second fundamental form is given by $2e^{-4v}$ and
then (\ref{gGauss}) implies that the induced metric by $\Phi^\vartheta $ is given by $e^{2u}|dz|^2$. In addition,
(\ref{CGauss}) says that the associated Jacobian to $\Phi^\vartheta $ is $C$.
Using a complex coordinate on $\Sigma$, we can write (see section 3.2) $\Phi^\vartheta  =  -2\sqrt 2 \,i e^{-2v} \, \Psi^\vartheta_z \wedge \Psi^\vartheta_{\bar z}$  and
$\widehat{\Phi^\vartheta}=  \sqrt{2} \, \Psi^\vartheta \wedge N^\vartheta $ and it is not difficult to get that 
the holomorphic 2-differential associated to $\Phi^\vartheta $ is $\Theta^\vartheta (z)=-i e^{i\vartheta}(dz)^2$. 
In addition, equation (\ref{sinhGordon}) means that the functions $u$ and $C$ satisfy the compatibility 
equations (\ref{GCeqns}). Therefore, our immersion $\Phi$ is congruent to the Gauss map $\Phi^{\theta+\pi/2}$ of $\Psi^{\theta+\pi/2}$.
\end{proof}

\subsection{A distinguished example}
We emphasize an interesting example of a Klein bottle studied in [EGJ] and [JNP] whose 
double cover is an $\s^1$-equivariant minimal torus in $\s^4$. 
Up to congruences, we are going to look at it as a minimal Lagrangian Klein bottle embedded in $\s^2 \times \s^2$.
We consider $\s^2=\{(x,z)\in\r\times\c\,/\,x^2+|z|^2=1\}$ and define 
\begin{equation}\label{defB}
{\bf B}=\{((x,z),(y,w))\in \s^2\times\s^2 \,/\,2x=y, \ \Re(\sqrt{z}w)=\Im(\sqrt{z}w)\},
\end{equation}
where $\sqrt{}$ stands for the main branch of the square root. 
We must point out that $|z|^2=1-y^2/4\geq 3/4$. 
We are able to give a conformal parametrization of the universal covering of 
$\bf B$ by means of elementary Jacobi elliptic functions (see [BF] for background):
\begin{equation}\label{confparamB}
\Phi_{\bf B}=(\phi,\psi):\r^2\rightarrow\s^2\times\s^2
\end{equation}
with
\begin{eqnarray*}
\phi(t,s)=\frac{1}{3\,\hbox{dn}(\sqrt{3}t)}\left(-\hbox{sn}(2\sqrt{3}t),
				i(2\,\hbox{sn}^2(\sqrt{3}t)-3)e^{\frac{4is}{\sqrt{3}}}\right ),\\
\psi(t,s)=\frac{1}{3\,\hbox{dn}(\sqrt{3}t)}\left(-2\,\hbox{sn}(2\sqrt{3}t),
				-i(4\,\hbox{sn}^2(\sqrt{3}t)-3)e^{\frac{-2is}{\sqrt{3}}}\right ),
\end{eqnarray*}
where sn, cn and dn stand for the sine amplitude, the cosine amplitude and the delta amplitude with modulus
$p=2\sqrt 2 /3$.

The corresponding group of transformations in $\r^2$ which defines ${\bf B}$ is generated by 
\[
(t,s)\mapsto\left(t+\frac{2\sqrt{3}K}{3},s\right), \quad (t,s)\mapsto\left(\frac{\sqrt{3}K}{3}-t,s+\frac{\sqrt{3}\pi}{2}\right),
\]
where $K$ is the complete elliptic integral of the first kind with modulus $p=2\sqrt 2 /3$.

We remark that $\Phi_{\bf B}$ is the Gauss map of the minimal immersion $\Psi:\r^2\rightarrow\s^3$ given by
\[
\Psi(t,s)=\left( \hbox{cn}(\sqrt{3}t)\,e^{i\sqrt{3}s},\hbox{sn}(\sqrt{3}t)\,e^{\frac{is}{\sqrt{3}}}\right).
\]
This is also a conformal parametrization of the universal covering of Lawson's $\tau_{3,1}$ torus described in [L]
by the orthogonal parametrization
\(
(u,v) \mapsto (\cos v \,e^{3iu}, \sin v \,e^{iu}).
\)
In addition, the Gauss map of Lawson's $\tau_{3,1}$ torus is a four-fold covering of the double cover torus of $\bf B$.
Following the notation of the proof of Theorem \ref{Th:min}, Lawson's $\tau_{3,1}$ torus corresponds to the solution 
$v=v(t)=\log (\sqrt 3 \,\hbox{dn}(\sqrt 3 t))$ of equation
(\ref{sinhGordon}) depending on only one variable, satisfying $\tanh v(0)=1/2, \, v'(0)=0$ and choosing $\vartheta = \pi/2$ (or $\theta =0$).

We point out that $\Phi_{\bf B} $ is a minimal immersion in 
$\s^4(\sqrt{2})$ (note that $\bf B$ lies in the hyperplane $2x=y$ of $\r^6$) and, making use of (\ref{confparamB}), 
it is not difficult to compute the area of $\bf B$ using the above data.  A straightforward computation leads to Area$({\bf B})=12 \pi E$, where $E$ is the complete elliptic integral of the second kind with modulus $p=2\sqrt 2 /3$.
Then Theorems 1.3.1 and 1.4.1 in [JNP] show that the first positive eigenvalue $\lambda_1$ of the Laplacian $\Delta$ (acting on functions) of the Klein bottle $\bf B$ is $\lambda_1({\bf B})=1$.

\section{Second variation of minimal Lagrangian surfaces}

 Let $\Phi:\Sigma\rightarrow\s^2\times\s^2$ be a minimal Lagrangian immersion of a compact surface $\Sigma$. 
We identify the sections on the normal bundle $\Gamma(T^{\perp}\Sigma)$ with the $1$-forms on $\Sigma$ by
\begin{eqnarray}\label{iden}
\Gamma(T^{\perp}\Sigma)&\equiv&\Omega^1(\Sigma)\\
\xi&\equiv& \alpha \nonumber
\end{eqnarray}
being $\alpha$ the $1$-form on $\Sigma$ defined by $\alpha(v)=\omega(\Phi_* v,\xi)$ for any $v$ tangent to $\Sigma$.
In this way, the Jacobi operator of the second variation of the area becomes in an intrinsic operator, which is given by (see [O])
\begin{eqnarray*}
L:\Omega^1(\Sigma)\rightarrow\Omega^1(\Sigma)\\
\alpha\mapsto\Delta\alpha+\alpha,
\end{eqnarray*}
where, in general,  $\Omega^p(\Sigma)$, $p=0,1,2$, is the space of $p$-forms on $\Sigma$ and $\Delta$ is the Laplacian of the induced metric, i.e. $\Delta=\delta d+d\delta$, where $\delta$ is the codifferential operator of the exterior differential $d$.

Hence, {\it the index of $\Phi$, that we will denote by $\hbox{Ind}\,(\Sigma)$, is the number of eigenvalues (counted with multiplicity) of $\Delta:\Omega^1(\Sigma) \rightarrow \Omega^1(\Sigma)$ less than $1$}.

In order to study the Jacobi operator, we  consider the Hodge decomposition
\[
\Omega^1(\Sigma)={H}(\Sigma)\oplus d\,\Omega^0(\Sigma)\oplus\delta\Omega^2(\Sigma),
\]  
which allows to write, in a unique way, any 1-form $\alpha$ as $\alpha=\alpha_0 + dg +\delta \beta$, being $\alpha_0$ a harmonic 1-form, $g$ a real function and $\beta$ a 2-form on $\Sigma$. The space of harmonic 1-forms, ${ H}(\Sigma)$, is the kernel of $\Delta$ and its dimension is the first Betti number $\beta_1(\Sigma)$ of $\Sigma$. 
As $\Delta$ commutes with $d$ and $\delta$, the positive eigenvalues of $\Delta:\Omega^1(\Sigma) \rightarrow \Omega^1(\Sigma)$ are the positive eigenvalues of $\Delta:\Omega^0(\Sigma)\rightarrow\Omega^0(\Sigma)$ joint to the positive eigenvalues of $\Delta:\Omega^2(\Sigma)\rightarrow\Omega^2(\Sigma)$. Therefore
\begin{equation}\label{Ind}
\hbox{Ind}\,(\Sigma)=\beta_1(\Sigma)\,+\,\hbox{Ind}_0(\Sigma)\,+\,\hbox{Ind}_1(\Sigma),
\end{equation}
where {\it $\hbox{Ind}_0(\Sigma)$ is the number of positive eigenvalues (counted with multiplicity) of $\Delta:\Omega^0(\Sigma)\rightarrow \Omega^0(\Sigma)$  less than $1$} and {\it $\hbox{Ind}_1(\Sigma)$ is the number of positive eigenvalues (counted with multiplicity) of $\Delta:\Omega^2(\Sigma)\rightarrow \Omega^2(\Sigma)$ less than $1$}.  

The variational vector fields of the Hamiltonian deformations of the Lagrangian surface $\Sigma$ are the normal components of the Hamiltonian vector fields on $\s^2\times\s^2$. 
If $F:\s^2\times\s^2\rightarrow\r$ is a smooth function and $X=J\bar\nabla F$ is its associated Hamiltonian vector field on $\s^2\times\s^2$, the $1$-form associated to the normal component of $X$, via the identification (\ref{iden}), is $d(F\circ\Phi)$. Thus our minimal Lagrangian compact surface $\Sigma$ is {\em Hamiltonian stable}, i.e. stable under Hamiltonian deformations, if the first positive eigenvalue of $\Delta$ acting on $\Omega^0(\Sigma)$ is at least $1$. But from (\ref{mcR6}) we have that $\Phi:\Sigma\rightarrow\s^5(\sqrt 2)$ is also a minimal immersion and so $\Delta\Phi+\Phi=0$, i.e. $1$ is an eingenvalue of $\Delta$. Hence {\it $\Sigma$ is Hamiltonian stable if the first positive eigenvalue of $\Delta$ acting on $\Omega^0(\Sigma)$ is $1$}.

Precisely the first eigenvalue of $\Delta$ (acting on functions) of the Lagrangian sphere ${\bf M}_0$ defined in (\ref{defM0}) is $1$ and it is clear that the same happens to the Lagrangian torus $\bf T$ defined in (\ref{defT}).
In section 4.4 we showed the same property for the Klein bottle $\bf B$ defined in (\ref{defB}). 
As a consequence, we have that ${\bf M}_0$, $\bf T$ and $\bf B$  are Hamiltonian stable minimal Lagrangian surfaces in $\s^2 \times \s^2$. 

If the compact surface $\Sigma$ is orientable, the star operator $\star:\Omega^0(\Sigma)\rightarrow\Omega^2(\Sigma)$ says us that the eigenvalues of $\Delta$ acting on $\Omega^0(\Sigma)$ or on $\Omega^2(\Sigma)$ are the same, and so $\hbox{Ind}_0(\Sigma)=\hbox{Ind}_1(\Sigma)$. Thus if $\Sigma$ is a minimal Lagrangian compact {\it orientable} surface of $\s^2\times\s^2$ with genus $g$, then 
\begin{equation}\label{Indor}
\hbox{Ind}\,(\Sigma)=2\,g\,+\,2\,\hbox{Ind}_0(\Sigma).
\end{equation}
Using (\ref{Indor}) we get that ${\bf M}_0$ is stable and $\hbox{Ind}({\bf T})=2$. 

In the following result we provide variational characterizations of the examples ${\bf M}_0$, $\bf T$ and $\bf B$ in this context.
\begin{theorem}\label{Th:indice}
Let $\Phi:\Sigma\rightarrow\s^2\times\s^2$ be a minimal Lagrangian immersion of a compact surface $\Sigma$. Then
\begin{enumerate}
\item If $\Sigma$ is stable, then $\Phi(\Sigma)$ is the totally geodesic Lagrangian sphere ${\bf M}_0$.
\item If $\Sigma$ is Hamiltonian stable and $\Sigma$ is orientable with genus $g\leq 2$, 
then $\Phi$ is an embedding and $\Phi(\Sigma)$ is either the totally geodesic sphere ${\bf M}_0$ 
or the totally geodesic torus $\bf T$.
\item If $\Sigma$ is a Hamiltonian stable Klein bottle, then $\Phi$ is an embedding and $\Phi(\Sigma)$ is the Klein bottle 
$\bf B$ described in section 4.4.
\item If $\Sigma$ is unstable, then $\hbox{Ind}(\Sigma)\geq 2$ and the equality holds if and only if $\Phi$ is an embedding and $\Phi(\Sigma)$ is the totally geodesic torus $\bf T$.
\end{enumerate}
\end{theorem}
\begin{remark}
{\rm Part (1) shows that the result proved in Corollary 5.2 in [MW] is the best possible. 
Also, since the totally geodesic sphere ${\bf M}_0$ is a complex surface with respect to the complex structure $J=(J,-J)$ on $\s^2\times\s^2$, we have that ${\bf M}_0$ is area minimizing in its homology class. In addition, in [IOS] it was proved that the totally geodesic torus $\bf T$ is area minimizing under Hamiltonian deformations of $\s^2\times\s^2$. }
\end{remark}
\begin{proof}
We start remenbering a result of Simon [S] which will be used in the proof of this Theorem. 

Let $\Psi:M\rightarrow\r^n$ be an immersion of a compact surface $M$ with mean curvature vector $\bar {H}$ and maximum multiplicity $\mu$, i.e. there exist $\mu$ points $\{p_1,\dots,p_{\mu}\}$ on $M$ such that $\Psi(p_i)=a$, for all $1\leq i\leq \mu$. Then
\[
\int_{M}|\bar{H}|^2dA\geq 4\pi\mu,
\]
and the equality holds if and only if  $\bar{H}$ is given on $\tilde{M}=M-\{p_1,\dots,p_{\mu}\}$ by $\bar{H}=\frac{-2(\Psi-a)^{\perp}}{|\Psi-a|^2}$, where $\perp$ stands for normal component. This condition about the mean curvature $\bar{H}$ means that $\frac{\Psi-a}{|\Psi-a|^2}:\tilde{M}\rightarrow\r^n$ is a minimal immersion.

In this setting, the minimal Lagrangian immersion $\Phi$ produces an immersion $\Phi:\Sigma\rightarrow\r^6$ which is minimal into $\s^5(\sqrt {2})$. From (\ref{mcR6}), in this case $\bar{H}=-\Phi/2$ and we obtain that
\[
\hbox{Area}\,(\Sigma)\geq 8\pi\mu,
\]
and the equality holds if and only if $\Phi=4\frac{(\Phi-a)^{\perp}}{|\Phi-a|^2}$, where $\mu$ is the maximum multiplicity of $\Phi$. Since now $a\in\s^5(\sqrt 2)$ and $\Phi$ is normal to the surface, the last equation becomes in $\langle\Phi,a\rangle\Phi=2a^{\perp}$. From here it is not difficult to conclude that $\Phi$ is the totally geodesic sphere ${\bf M}_0$. As a summary, we have obtained that
 \begin{equation}\label{EstArea}
\hbox{Area}\,(\Sigma)\geq 8\pi\mu
\end{equation}
and {\em the equality holds if and only if $\Phi $ is an embedding and $\Phi(\Sigma)$ is congruent to ${\bf M}_0$, whose area is $8\pi$.}

Now we can prove (1). From (\ref{Ind}), if $\Sigma$ is stable then $\beta_1(\Sigma)=0$, which implies, taking into account Proposition \ref{Prop:min+cpt+or},(2), that $\Phi(\Sigma)$ is congruent to ${\bf M}_0$.

 To prove (2), we only have to consider the cases $g=1,2$. If the genus of $\Sigma$ is $1$, i.e. $\Sigma$ is a Hamiltonian stable torus, then the first eigenvalue of the Laplacian acting on functions is $1$. So we have a minimal immersion $\Phi:\Sigma\rightarrow\s^5(\sqrt{2})$ of a torus $\Sigma$ where $1$ is the first eigenvalue of the Laplacian. A result of El Soufi and Ilias [EI] says us that $\Sigma$ is flat and then Theorem \ref{Th:min+K} shows that it must be the totally geodesic torus $\bf T$. 

If the genus of $\Sigma$ is $2$ we use a known argument. From the Brill-Noether theory, we can get a non constant meromorphic map $\varphi:\Sigma\rightarrow\s^2$ of degree $d\leq 2$. Then there exists a Moebius transformation $F:\s^2\rightarrow\s^2$ such that $\int_{\Sigma}(F\circ\varphi)=0$, and using that the first positive eigenvalue of $\Delta$ is $1$, we have  
\[
\int_{\Sigma}|\nabla(F\circ\varphi)|^2\geq \int_{\Sigma}|F\circ\varphi|^2=\,\hbox{Area}(\Sigma).
\]
But $\int_{\Sigma}|\nabla(F\circ\varphi)|^2=8\pi\,\hbox{degree}(F\circ\varphi)=8\pi\,\hbox{degree}(\varphi)\leq 16\pi$. Hence we obtain that $\hbox{Area}(\Sigma)\leq 16\pi$. Taking into account (\ref{EstArea}) we finally get that
\(
\mu\leq 2.
\)
But Proposition \ref{Prop:or+cpt+emb} implies that $\mu\geq 2$, obtaining the equality in the last inequality. This gives a contradiction, because the equality can be only attained by the totally geodesic ${\bf M}_0$. This finishes the proof of (2).
 
Suppose now that $\Sigma$ is a Hamiltonian stable Klein bottle. Then we have a minimal immersion $\Phi:\Sigma\rightarrow\s^5(\sqrt{2})$ of a Klein bottle $\Sigma$ such that $1$ is the first positive eigenvalue of $\Delta$. From Theorem 1.2 in [EGJ] we deduce that our immersion is an embedding and the surface is the Klein bottle $\bf B$. This proves (3).

Finally we prove (4). If $\Sigma$ is unstable and orientable, using part (1) the genus $g$ of $\Sigma $ satisfies that 
$g\geq 1$ and, from (\ref{Indor}), $\hbox{Ind}\,(\Sigma)\geq 2$ and the equality holds if and only if $\Sigma$ is a Hamiltonian stable torus, which implies that it is the totally geodesic torus $\bf T$ using part (2).

If $\Sigma$ is non orientable, let $\pi:\tilde{\Sigma}\rightarrow\Sigma$ the $2:1$ orientable Riemannian covering and $\tau:\tilde{\Sigma}\rightarrow\tilde{\Sigma}$ the change of sheet involution. The spaces of forms on $\tilde{\Sigma}$ can be decomposed in the following way:
\[
\Omega^i(\tilde{\Sigma})=\Omega^i_+(\tilde{\Sigma})\oplus\Omega^i_-(\tilde{\Sigma}),\quad i=0,1,2,
\]
where
\[
\Omega^i_{\pm}(\tilde{\Sigma})=\{\alpha\in\Omega^i(\tilde{\Sigma})\,/\,\tau^*\alpha=\pm \alpha\}.
\]
As $\pi\circ\tau=\pi$, the map $\alpha\in\Omega^i(\Sigma)\mapsto\pi^*\alpha\in\Omega^i(\tilde{\Sigma})$ allows to identify $\Omega^i(\Sigma)\equiv\Omega^i_+(\tilde{\Sigma})$, $i=0,1,2$. In addition, as $\Sigma$ is non orientable, $\star\tau^*=-\tau^*\star$, and so $\star$ identifies $\Omega^0_-(\tilde{\Sigma})\equiv\Omega^2_+(\tilde{\Sigma})$. Hence we have the identification 
\begin{eqnarray*}
\Omega^2(\Sigma)&\equiv&\Omega^0_-(\tilde{\Sigma})\\
\beta&\equiv& f  
\end{eqnarray*}
where $\pi^*\beta=f\omega_{\tilde{\Sigma}}$, being $\omega_{\tilde{\Sigma}}$ the area $2$-form on $\tilde{\Sigma}$.
Since $\Sigma$ is non orientable, the eigenvalues of $\Delta:\Omega^2(\Sigma)\rightarrow\Omega^2(\Sigma)$ are positive, and so, taking into account the above considerations, {\em $\hbox{Ind}_1(\Sigma)$ is the number of eigenvalues (counted with multiplicity) of $\Delta:\Omega^0_-(\tilde{\Sigma})\rightarrow\Omega^0_-(\tilde{\Sigma})$ less than $1$}. Also, as {\em $\hbox{Ind}_0(\Sigma)$ is the number of positive eigenvalues (counted with multiplicity) of $\Delta:\Omega^0_+(\tilde{\Sigma})\rightarrow\Omega^0_+(\tilde{\Sigma})$ less than $1$}, we obtain that
\begin{equation}\label{FlaInd}
\hbox{Ind}_0(\Sigma)+\hbox{Ind}_1(\Sigma)=\hbox{Ind}_0(\tilde{\Sigma}),
\end{equation}
corresponding to the minimal Lagrangian immersion $\Phi\circ\pi:\tilde{\Sigma}\rightarrow\s^2\times\s^2$.

We can consider that $\chi(\Sigma)\leq 0$ because there are not minimal Lagrangian real projective planes in $\s^2\times\s^2$ according to Proposition \ref{Prop:min+cpt+or},(2). 
If $\chi(\Sigma)\leq -2$, then $\beta_1(\Sigma)\geq 3$ and (\ref{Ind}) says that $\hbox{Ind}\,(\Sigma)\geq 3$. If $\chi(\Sigma)=-1$, then $\beta_1(\Sigma)=2$ and 
(\ref{Ind}) and (\ref{FlaInd}) imply that
\[
\hbox{Ind}\,(\Sigma)=2+\hbox{Ind}_0(\Sigma)+\hbox{Ind}_1(\Sigma)=2+\hbox{Ind}_0(\tilde{\Sigma}).
\]
But $\tilde{\Sigma}$ is an oriented compact surface of genus 2 and part (2) leads to 
$\hbox{Ind}_0(\tilde{\Sigma})\geq 1$. This implies that $\hbox{Ind}\,(\Sigma)\geq 3$.

Finally, if $\Sigma$ is a Klein bottle, then $\beta_1(\Sigma)=1$. 
If $\hbox{Ind}_1(\Sigma)=0$, using the above description of  $\hbox{Ind}_1$,
the first eigenvalue $\lambda_1$ of $\Delta:\Omega^0_-(\tilde{\Sigma})\rightarrow\Omega^0_-(\tilde{\Sigma})$ satisfies $\lambda_1\geq 1$. Hence 
\[
\int_{\tilde{\Sigma}}|\nabla f|^2\geq \int_{\tilde{\Sigma}}f^2,\quad \forall f\in C^{\infty}(\tilde{\Sigma})\quad\hbox{such that}\quad f\circ\tau=-f.
\]
From Theorem 1 in [RS], we can get a non constant meromorphic map $\varphi:\tilde{\Sigma}\rightarrow\s^2$ satisfying $\varphi\circ\tau=-\varphi$ of degree $d\leq 2$. Thus we obtain 
\[
\int_{\tilde{\Sigma}}|\nabla \varphi|^2\geq \int_{\tilde{\Sigma}}|\varphi|^2=\hbox{Area}(\tilde{\Sigma}).
\]
But $\int_{\tilde{\Sigma}}|\nabla \phi|^2=8\pi\,\hbox{degree}(\phi)\leq 16\pi$. So we get that $\hbox{Area}(\tilde{\Sigma})\leq 16\pi$. Since $\hbox{Area}(\tilde{\Sigma})=2\hbox{Area}(\Sigma)$, from (\ref{EstArea}) we have that $\mu\leq 1$. Hence $\mu=1$ and the equality in (\ref{EstArea}) holds, which is impossible because $\tilde{\Sigma}$ would be the totally geodesic sphere ${\bf M}_0$. Therefore $\hbox{Ind}_1(\Sigma)\geq 1$ for any minimal Lagrangian Klein bottle $\Sigma$ of $\s^2\times\s^2$. Part (3) gives that $\hbox{Ind}_0(\Sigma)\geq 1$ for any minimal Lagrangian Klein bottle of $\s^2\times\s^2$ except for {\bf B}. In this way, we obtain that any minimal Lagrangian Klein bottle $\Sigma$ different from {\bf B} satisfies $\hbox{Ind}(\Sigma)\geq 3$. 

To finish the proof we must check that $\hbox{Ind}_1({\bf B})\geq 2$. In fact, following section 4.4, we consider the functions
$f,g:\tilde{{\bf B}}\rightarrow\r$ defined by $f(\pi(t,s))=\cos(2s / \sqrt 3)$, $g(\pi(t,s))=\sin (2s/\sqrt 3)$, 
where $\pi:\r^2\rightarrow\tilde{{\bf B}}$ is the projection and $\tilde{{\bf B}}$ is the $2:1$ covering torus of {\bf B}. As the involution $\tau:\tilde{{\bf B}}\rightarrow\tilde{{\bf B}}$ is induced by 
$(t,s)\mapsto(\sqrt{3}K /3-t,s+\sqrt{3}\pi /2)$, it is clear that $f\circ\tau=-f$ and $g\circ\tau=-g$. Following the proof of Theorem \ref{Th:min}, the induced metric on $\r^2$ is given by $e^{2u(t)}=4\cosh\log(\sqrt{3}\,\hbox{dn}(\sqrt{3}t))$ and so $e^{2u(t)}\geq 4$.
Then, for any real numbers $a$ and $b$ we have that
\[
\begin{array}{c}
\Delta (af+bg)(\pi(t,s))=e^{-2u(t)}\frac{d^2}{ds^2}\left(a\cos\frac{2s}{\sqrt{3}}+b\sin\frac{2s}{\sqrt{3}}\right)\\
=-\frac{4}{3}e^{-2u(t)}(af+bg)(\pi(t,s)).
\end{array}
\]
But using that $e^{-2u(t)}\leq 1/4$, it follows that
\[
-(af+bg)\Delta(af+bg)\leq\frac{1}{3}\left( af+bg \right)^2,
\]
which implies that 
\[
-\int_{{\bf B}}(af+bg)L(af+bg)\,dA\leq-\frac{2}{3}\int_{{\bf B}}\left( af+bg \right)^2\,dA.
\]
In conclusion, we have shown that there exists a $2$-dimensional subspace of $\Omega^0_{-}(\tilde{{\bf B}})$ whereon the quadratic form associated to the Jacobi operator $\Delta+1$ of $\tilde{{\bf B}}$ is negative definite. So $\hbox{Ind}_1({\bf B})\geq 2$ and this finishes the proof.
\end{proof}
\vspace{0.3cm}

In the last result we compute the index of the Gauss map of a compact orientable minimal surface of $\s^3$ in terms of the index of itself. 
\begin{proposition}
Let $\Psi:\Sigma\rightarrow\s^3\subset\r^4$ be a minimal immersion of an orientable compact surface $\Sigma$ and $\Phi:\Sigma\rightarrow\s^2_+\times\s^2_-$ its Gauss map. Then
\[
\hbox{Ind}_0\,(\Phi)=\hbox{Ind}\,(\Psi)-1.
\]
Moreover, if the genus of the surface $g\geq 1$ then $\hbox{Ind}\,_0(\Phi)\geq 4$ and $\hbox{Ind}\,(\Phi)\geq 10$, 
and the equality holds in some of the equalities if and only if $\Phi$ is a two fold covering of the totally geodesic Lagrangian torus $\bf T$ (see Remark \ref{Re:Clifford}).
\end{proposition}
\begin{proof}
The Jacobi operator of the second variation of $\Psi$ is given by $L=\hat{\Delta}+|\hat{\sigma}|^2+2$ (see [U]), where $\hat{\Delta}$ is the Laplacian of the induced metric $\hat{g}$. But from (\ref{gGauss}) the induced metric $g$ by the immersion $\Phi$ is conformal to $\hat{g}$ with $g=(2+|\hat{\sigma}|^2)\hat{g}$. Hence the quadratic form $\hat{Q}$ associated to $\hat{L}$ acting on functions of $\Sigma $ verifies
\[
\begin{array}{c}
\hat{Q}(u,u)=-\int_{\Sigma}\{u\hat{\Delta}u+(|\hat{\sigma}|^2+2)u^2\}dv_{\hat{g}}=\\
=-\int_{\Sigma}(u\Delta u+u^2)(|\hat{\sigma}|^2+2)dv_{\hat{g}}=-\int_{\Sigma}(u\Delta u+u^2)dv_g=Q(u,u),
\end{array}
\]
where $Q$ is the quadratic form associated to the operator $\Delta+1$. So $\hbox{Ind}_0\,(\Phi)=\hbox{Ind}\,(\Psi)-1$, because to compute $\hbox{Ind}_0\,(\Phi)$ we only consider positive eigenvalues of $\Delta$. Now, we use the main result proved in [U]: We have that if $g\geq 1$, then $\hbox{Ind}\,(\Psi)\geq 5$ and the equality holds if and only if $\Psi$ is the Clifford torus. From (\ref{Indor}) and the fact that the Gauss map of the Clifford torus is the two fold covering of $\bf T$ (see Remark \ref{Re:Clifford}), we finish the proof.
\end{proof}

\vspace{1cm}

\end{document}